# The Theory of Fallible Probability

## and

## The Dynamics of Degrees of Belief

Amos Nathan[1]

amosnathan1@gmail.com

**ABSTRACT**

This monograph is an account of the theory of fallible probability and of the dynamics of degrees of belief. It discusses the first order 'subjective' theory in which first order degrees of belief are expressed by 'subjective' probabilities and are updated by conditionalization (Bayes, 1764; Ramsey, 1926), gives an improved exposition of the greater part of the author's theory of 'Probability Dynamics' (Nathan, 2006) which should replace the so-called 'Probability Kinematics' (Jeffrey, 1965), resolves the problem of 'New Explanation of Old Evidence' (Jeffrey, 1995), provides a Theory of Confirmation, and refutes the 'Principle of Reflection' (Van Fraassen, 1984).

# Table of Contents

**SURVEY**

Part I analyzes the conventional first order method of updating by conditionalization which has met a spate of objections. It proves its soundness and thereby rebuts all and any counterexamples, including the 'Problem of Old Evidence' and some supposed difficulties with degrees of confirmation.

Part II proves that there exist no non-vanishing second order probabilities of first order epistemic degrees of belief, shows that inconsistencies of higher order probability are not eliminated by constraints that have been proposed for this purpose.

Part III develops the theory of updating fallible probability, as first given in the author's 'Probability Dynamics' (PD)[4], which quantifies relative degrees of relevance of fallible first order probability distributions by means of an associated 'credence' thereby allowing the formation of mergers of probabilistic propositions such as '$P(A)=p$' and '$P(A)=q$'. PD should replace the theory of 'Probability Kinematics' (Jeffrey, 1965) which is shown to be inconsistent. PD is 'Straight' (SPD) or 'Offsetting' (OPD). SPD produces mutually reinforcing 'straight' mergers, while mergers of OPD set off conflicting evidences against each other according to their relative credence and their correlation. Also discussed is a process of 'normalization' that provides an interpretation of sets of binary evidences of non-uniform credence as equi-credible distributions; a 'cross-credence' that determines the credence of interacting evidences; and 'indirect' revisions of evidence for one distribution by evidence for another. This Part includes some new material and corrects an error in the author's earlier treatment of indirect updating.

Part IV points out that the 'Problem of New Explanation of Old Evidence' cannot be solved by first order methods, resolves it, and disproves the earlier results of Jeffrey (1995) and his followers.

Part V is a PD based theory of confirmation (PDCT).

Part VI refutes the so-called Principle of Reflection.

Part I, which discusses first order subjective probability, and Part II are no essential prerequisites for Parts III to V. Part V depends just on Part III.

**PART I      The First Order Theory**

**1.1      Conditionalization**

In the theory of 'subjective', 'epistemic', or 'propositional', probability, agents are assumed to be rational in the sense of maximizing their personal gain, and $P(A)$, $P(B)$ etc. are an agent's degrees of belief in propositions A, B, etc. The theory accepts the  pragmatic definition of degrees of belief as agents' betting ratios in bets on propositions (Bayes, 1763; Ramsey 1926; De Finetti 1970, 1972; Savage 1972; Skyrms 1980, 1986, 1987). Such degrees of belief satisfy the Kolmogorov axioms of probability (Gnedenko, 1989).

A proposition is *given* if it is believed to be true and *accepted* if it is newly given.  Probability $P_0(A)$ is $P(A)$ at time $t_0$, $P_1(A)$ is $P(A)$ at $t_1$, 'prior' $E$ is the collection of probabilities $E=\{P(A),P(AB),P(B)\}$,

$$E_0=\{P_0(A), P_0(AB), P_0(B)\} \tag{1.1}$$

is $E$ at $t_0$, etc., and, contrary to conventional usage, we *define* $P(A|B)$ by $P(A|B)P(B)=P(AB)$, $P_0(A|B)$ by $P_0(A|B)P_0(B)=P_0(AB)$, etc. [2]

Casting myself as a rational agent, my probability $P(A)$ of A is my equitable betting ratio for a bet on A, i.e., my equitable price for a bet on A that yields 1 if A comes true and nothing otherwise. My probability $P_{01}(A/B)$ of A conditional on B is my equitable betting ratio at time $t_0$ (first subscript of $P_{01}$) for a bet on A that is on if B comes true at some later time $t_1$ (second subscript of $P_{01}$) and is called off otherwise [Ramsey, 1926]. Thus

$$P_{01}(A/B) = P_0(\overline{B})P_{01}(A/B) + P_0(AB),$$

i.e.,

$$\text{given } E_0, \quad P_{01}(A/B) = P_0(A|B) \text{ if } P_0(B) > 0, \tag{1.2}$$

which satisfies the necessary condition $P_{01}(A/B) + P_{01}(\overline{A}/B) = 1$. $t_1$ must not be predetermined, for if it were, then $P_{01}(A/B)$ would vanish identically.[3]

Note that $P_{01}(A/B)$ is my probability of A conditional on B and neither my probability of A given B nor my probability of A given my belief that B will come true.

Let $E_{10}$ at $t_0$ (second subscript of $E_{10}$) be the prior I believe at $t_1$ to have been my last update of $E$, set $P_{10}(A|B)=P_{10}(AB)/P_{10}(B)$, and let $P_{11}(A/B)$ be my probability of A at $t_1$ (first subscript of $P_{11}$) conditional on B at $t_1$ (second subscript of $P_{11}$). With these definitions it follows from equ. (2) that

$$\text{given } E_{10}, \quad P_{11}(A/B) = P_{10}(A|B) \text{ if } P_{10}(B) > 0. \tag{1.3}$$

If, in addition, $E_{10}=E_0$, then $P_{10}(A|B)=P_0(A|B)$ and

$$\text{given } E_{10}, \quad P_{11}(A/B)=P_0(A|B) \text{ if } P_0(B)>0. \tag{1.4}$$

If I accept B at $t_1>t_0$ then $P_{11}(A/B)$ is my update $P_1(A)$ at $t_1$ of $P_0(A)$ at $t_0$. Therefore, if I accept B at $t_1>t_0$ and $E_{10}=E_0$, then

$$\text{given } E_0, \quad P_1(A)=P_0(A|B) \text{ if } P_0(B)>0. \tag{1.5}$$

Equ. (5) has been challenged by a spate of putative counterexamples. If our arguments are sound, equ. (4) admits no exceptions and any conceivable counterexample to equ. (5) is reduced to the claim that $E_{10}$ is not necessarily equal to $E_0$.

The theory under consideration is formal. It postulates that agents are rational and is not otherwise concerned with attributes of agents. Once degrees of belief have been assigned to propositions any further development is purely syntactic.

New evidence can change my beliefs about the past but cannot revise the past. Therefore I know that my acceptance of B at $t_1>t_0$ cannot revise my belief $E_0$ at $t_0$ and will not revise my belief $E_{10}$ at $t_1$ of what $E_0$ has been at $t_0$. Hence $E_{10}=E_0$ and first order indirect updating always proceeds by equ. (5).

This result allows us to take $P_0(A|B)$ for $P_{01}(A/B)$.

## 1.2 The Problem of Old Evidence

The so-called 'Problem of Old Evidence' (Glymour 1980) is an abundantly discussed counterexample to equ. (5). (e.g., by Gaifman, 1985; Earman, 1992; Christensen, 1999; Barnes, 1999). It has led Glymour not just to "pessimistic conclusions for the prospects of Bayesian confirmation theory" but to the conclusion that "the argument, if sound, is destructive of the whole theory of subjective probability." Earman's discussion includes a survey of several ineffectual attempts to resolve the problem and he is likewise skeptical, and so are Bogen (1995)[4] and many others. The problem is distinct from the Problem of New Explanation of Old Evidence (Jeffrey, 1995) which is discussed in Part IV.

This is the Problem: Suppose that B is first accepted at $t_0$, i.e., suppose that $P_0(B)=1$. Then B is old evidence at $t_1>t_0$, and by equ. (5) $P_1(A)=P_0(A|B)=P_0(A)$, which has been said to prove that old evidence cannot update P(A), whereas, surely, it should be able to update P(A) no less than new evidence.

Let us mention one more of its many expositions. According to Christensen there is a "fatal flaw" in the theory, since 'P(A|B)=P(A) if P(B)=1' means that "if one is absolutely certain about one's evidence, it no longer is evidence." He points out that such a B is a tautology and says that "we should never believe a nontautology with probability 1 anyway. But this response does not remove the problem's sting. For it remains true that as P(B) approaches 1, the degree to which B can confirm anything becomes vanishingly small."

To resolve the problem we metaphrase Christensen: -- As $P_0(B)$ increases, the degree to which a B accepted at $t_1$ updates P(B), decreases, and if $P_0(B)=1$ then B accepted at $t_1$ merely repeats what has been accepted at $t_0$ and is redundant. In other words, $P_0(B)=1$ pre-empts at $t_0$ the repeated acceptance of B at $t_1$. Moreover, $P_0(B)=1$ entails that A and B are statistically independent, and a probability of A is not to be updated by a B of which it is statistically independent. The supposed problem confuses 'P(B)=1' accepted at $t_1$, with $P_0(B)$ given at $t_0$ as a component of prior $E_0$.

Many discussions of the formation of $P_0(A|B)$ do not mention $E_0$ at all. In many others 'background knowledge' is said to be required, and this background knowledge, say C, is taken into account by writing $P_0(A|BC)$ for $P_0(A|B)$. This is no adequate improvement. In the first place, since C does not stand for just any background knowledge but is $E_0$, and, more significantly, B and C (=$E_0$) play different roles in $P_0(A|BC)$. While B is a *conditioning* clause, $E_0$ is *given*.

Far from being "destructive of the whole theory of subjective probability", 'old evidence' is treated by the theory precisely as it should.

If $P_0(B)=1$ and a last earlier prior $E^*$ with $0<P^*(B)<1$ has been given before, then component $P_0(A)$ of $E_0$ will have been produced by conditionalization on B at $t_0$ given $E^*$.

Updating by fallible $P_1(B)$ requires a theory of higher order.

## PART II      PHO:  Higher Order Probability

Degrees of belief that do not satisfy the sum and product rules of the probability calculus are inconsistent (Cox, 1946), therefore the only feasible higher order theories of degrees of belief may seem to be theories of higher order probability (PHO).  Several attempts have been made to provide them (e.g., Soshichi, 1973; Skyrms, 1980, 1986; Gaifman, 1985; Goodman and Nguyen, 1999). We prove [5]

1) that there exist no positive second order probabilities of epistemic first order probabilities;
2) that PHO produces mergers that are not associative and commutative as they would have to be in order to allow a clear semantic interpretation;
3) that inconsistencies of PHO are not eliminated by constraints which have been suggested for this purpose;
4) that PHO cannot properly account for counter-evidence.

Let P(A) be $p \in (0,1)$ be predetermined [6], and let an accepted evidence not preclude the existence of an $\varepsilon$-neighbourhood $D_\varepsilon \subset (0,1)$ of $p$ in which there exists an $\eta>0$ such that $P[P(A)=q]>\eta$ for every $q \in D_\varepsilon$. This is the case for all epistemic propositions A.  Since 'P(A)=$r$' and 'P(A)=$s$' are mutually exclusive unless $r=s$, the probability that P(A) is in $D_\varepsilon$ is the sum of all second order probabilities $P[P(A)=q]$ with $q \in D_\varepsilon$.  This sum is infinite and not $\leq 1$ as required by a probability.  Therefore $P[P(A)=p] \equiv 0$ if $p \in (0,1)$. This shows that there exist no probabilities of epistemic propositions but, more generally, that there exist no such absolute degrees of belief, but it does not exclude the existence of *relative* degrees of belief of probabilistic propositions.

The probability of the probability later tonight of rain tomorrow is always nil.[7] This is not true of the probability of the *prediction* later tonight of rain tomorrow, whose range is the alternative 'rain or no rain'.

We prove that Mergers by PHO are not associative and not commutative.

*Proof*

Let $a_i>0$ be the second order probability of $e_i$='P(A)=$p_i$', let $(a;e)$ be the PHO-merger of $(a_1;e_1),\ldots,(a_n;e_n)$ with $e$='P(A)=$p$', $e_i$='P(A)=$p_i$' and $a$ and $a_i$ the probabilities of $e$ and $e_i$, respectively.  If $(a;e)$ is associative and commutative then there exists a function $f(.)$ such that $f(x)>0$ and $a = \Sigma_1^n f(a_i)$. Since $\Sigma_1^n f(a_i)$ exceeds 1 for sufficiently great $n$ this is not generally possible and proves our assertion.

Let B='$P_1(A)=r$', $P_0(A|B)=s$, and $P_0(B) > 0$.   Then $r=0$ entails '$P_0(A|B)=P_0(A|\overline{A})=s=0$' and $r=1$ entails '$P_0(A|B)=P_0(A|A)=s=1$'.  In these two cases there holds the identity $P_0(A|P_1(A)=r)=r$.

With B='$P_1(A)=r$' as before, on the assumption that $P_0(A|B)$ exists and $0<r<1$,
  a) let $0<P_0(B)<1$ and let B be accepted at $t_1$. Then $P_1(A)=r>0$ and by equ. (1.5) $P_1(A)=P_0(A|B)=P_0(AB)/P_0(B)$.         Since $r<1$, $AB=(P(A)=1) \wedge (P_1(A)=r)=\varnothing$ and $P_1(B)=P_0(A|B)=0$. It follows that $P_0[A|P_1(A)=r]$ does not exist, for if it did, it had to

be simultaneously 0 and $r>0$.

b) let $P_0(B)=1$. Then $P_1(A)=r$ and $P_1(A)=P_0(A|B)=P_0(A)$, independent of B and therewith independent of $r$, in contradiction to '$P_1(A)=r$'. Therefore our assumptions are false and $P_0[A|P_1(A)=r]$ does not exist.

This result is brought about by the existence of a feedback loop in '$P_0(A|P_1(A)=r)=s$' and the binary nature of the first order '$P_1(A)=r$', in which '$P_0$' stands for '$P_{01}$'. Indeed, $s$ is $P_{01}(A)$, and this reflects back on and is equal to conditioning phrase '$P_1(A)=r$', which is consistent if and only if $r=s$.

In any discourse that includes no probabilistic feedback loops, second order *syntactic* probabilities are perfectly in order.

For every evidence $(a;e)$ there should exist a 'counter-evidence' $(-a;e)$ such that the two do cancel out, but in a theory of higher order degrees of belief it is not sufficient just to tack $(-a;e)$ on to $(a;e)$, instead, a continuously diminishing '$a$' should bring about a gradual transition from $(a;e)$ with a>0, via $(0;e)$ to $(-a;e)$, which is not possible if $a$ is a probability.

## Part III  Probability Dynamics (PD)

### 3.1  Introduction

We can simultaneously believe in some degree in propositions such as '$P(A)=p$' and '$P(A)=q$' even if $p$ and $q$ are unequal, but in a first order theory these propositions are not compatible unless they are the same. That there can be no absolute second order degrees of belief in such propositions does not preclude the existence of relative second order degrees of belief, and that there exist no higher order epistemic probabilities does not preclude the feasibility of a higher order theory of degrees of belief. Such a theory should be able to express degrees of fallibility of first order degrees of belief, allow their revision, and show how to produce mergers of distributions such as $(p_1,...,p_m)$ and $(q_1,...,q_m)$ on a common partition of a probability space.

The theory of Probability Dynamics[4] accepts the first order theory of Part I. It introduces a 'credence' $\kappa$ of probability distributions as a second order degree of belief. $\kappa$ is a relative measure of the credibility, reliability, or relevance of the evidence for a distribution. It ranges from $-\infty$ to $\infty$ and is relative to an arbitrary reference that remains fixed throughout any one discourse. The credibility of a first witness is $n$ times that of a second if his/her evidences E are equivalent to the same E independently given in evidence by $n$ second witnesses and, If E='$P(A)=p$' given in evidence by a first witness has credence $\kappa$ then the same E given by a second witness whose credibility is $n$ times that of the first has credence $n\kappa$. A distribution whose credence is nil is irrelevant, increasing credence expresses growing confidence in it, and negative credence designates it as counter-evidence.

Evidences are 'evidentially independent', or just 'independent', if they are independently given by mutually independent witnesses or report assessments of outcomes of pair-wise mutually independent sets of trials. The credence of $n$ identical evidentially independent evidences is proportional to $n$.

An objective credence is readily defined in a frequentist formulation of probability. Consider a coin-tossing game with outcomes H for heads and T for tails, and let probability P(H) be the probability of heads. Given a stationary time series of $n$ mutually independent trials, the frequentist theory of probability takes for P(H) the relative frequency $f_n$(H) of heads, with increasing confidence in the representation of P(H) by $f_n$(H) the greater the number of trials. In this representation 'straight' PD (SPD) makes $\kappa$ proportional to $n$.

On the other hand, if the outcome of a single toss in the coin-tossing example is fallibly reported by two mutually independent observers whose evidences are equally credible and one says heads while the other says tails, or if their testimonies regarding a given proposition are, respectively, 'true' and 'false', then the two should cancel out. These are extreme examples of 'offsetting' probability dynamics (OPD).

### 3.2 Elements of Probability Dynamics

PD operates on probability space $\{\Omega, \Sigma, P\}$, with $\Omega$ a finite set of discrete elementary events and P a function which assigns distributions to field $\Sigma$ of $\Omega$ that satisfy the laws of probability (Kolmogorov axioms: Gnedenko, 1989). PD augments $\{\Omega, \Sigma, P\}$ by assigning a 'credence' to distributions on $\Sigma$. Distributions confer their credence to all functions of their elements. PD provides a process of 'normalization' which interprets sets of probabilities of non-uniform credence as an equi-credible distribution.

'$\alpha$–Evidence' $\alpha=[\kappa;d;U_m]$ has credence $\kappa$ and distribution $d=(p_1,\ldots,p_m)$ on partition $U_m$ of $\Sigma$. Whenever adequate, $\alpha$ can be abbreviated by $\alpha=[\kappa;d]$ and $[\kappa;p,1-p]$ by $[\kappa;p]$. Evidences are 'evidentially independent' if they are independently provided in equivalent circumstances by pair-wise independent witnesses, or report observations of pair-wise mutually independent sets of trials.

PD is 'straight' ('SPD', or just 'PD') or offsetting ('OPD'). SPD merges pair-wise evidentially independent $\alpha$–evidences $\alpha_1=[\kappa_1;d_1]$, $\alpha_2=[\kappa_2;d_2]$ etc. on a common partition of $\Sigma$, by merging corresponding elements into 'straight' mergers, and OPD merges them into 'offsetting' mergers. SPD-mergers of evidences of positive credence are mutually reinforcing, associative and commutative. OPD-mergers set off conflicting evidences against each other. Frequentist SPD applies to outcomes of mutually independent trial runs and frequentist OPD applies to evidentially independent assessments of outcomes of a single trial run.

$\alpha_1 \oplus \ldots \oplus \alpha_n$, or $\{\alpha_i|_1^n\}$, is the straight merger of $\alpha_1,\ldots,\alpha_n$ and $\alpha_1 \Diamond \ldots \Diamond \alpha_n$ is their offsetting merger. With $(q_j, 1-q_j)$ a distribution on $\{A_j, \overline{A}_j\}$ for $j=1,\ldots,m$, let us define the set $\beta=\{[\kappa_1;q_1],\ldots,[\kappa_m;q_m]\}$. Here $(q_1,\ldots,q_m)$ is an $\alpha$–distribution on some partition $U_m=\{A_1,\ldots,A_m\}$ if and only if $\kappa_1 = \ldots = \kappa_m$ and $\Sigma_1^m q_j = 1$, but a process of 'normalization' interprets any set $\{[\kappa_1;q_1],\ldots,[\kappa_m;q_m]\}$ with non-negative $\kappa_j$ as an $\alpha$–evidence.

### 3.3 Straight Probability Dynamics (SPD)

#### 3.3.1 *Straight Mergers*

In SPD we prove the rules of merging binary $\alpha$-evidences of rational-valued non-negative credence of non-vanishing sum. Their extension to real valued credence presents no difficulty and will be omitted.

In a frequentist representation, on partition $U_2$, in a run of $n_i$ trials with $m_i$ successes and $n_i - m_i$ failures, let us set $\kappa_i = n_i$ and $p_i = m_i/n_i$. Then $\kappa = \Sigma n_i$, $p = \Sigma m_i / \Sigma n_i$, and $\{[\kappa_i; p_i]_1^n\} = [\kappa; p]$ with

$$\kappa = \Sigma_1^n \kappa_i; \qquad p = \Sigma_1^n \kappa_i p_i / \Sigma_1^n \kappa_i. \tag{3.1}$$

Alternatively, let us postulate that the straight merger of $n$ evidentially independent equi-credible evidences $[\kappa; p_1], \ldots, [\kappa; p_n]$ has credence $n\kappa$ and a probability equal to the average of the merged probabilities, i.e.,

$$\{[\kappa; p_i]_1^n\} = [n\kappa; \tfrac{1}{n} \Sigma_1^n p_i]. \tag{3.2}$$

Accordingly

$$\{[\kappa; p]_1^n\} = [n\kappa; p], \tag{3.3}$$

so that $n$ evidentially independent identical evidences $[\kappa; p]$ are equivalent to a single evidence of probability $p$ and credence $n\kappa$. We prove that equ. (2) entails equs. (1):

*Proof*

Let $\{[\kappa; d]_\oplus^n\}$ denote the $n$-fold merger of $[\kappa; d]$ with itself and let $\kappa_i$ be $r_i/s_i$ with integer $r_i$ and $s_i$, then by (2)

$$\{[\kappa_i; p_i]_1^n\} = \{[1/\Pi_1^n s_j; p_1]^{r_1 \Pi_{j\neq 1}^n s_j}\}_\oplus \oplus \ldots \oplus \{[1/\Pi_1^n s_j; p_n]^{r_1 \Pi_{j\neq n}^n s_j}\}_\oplus$$

$$= \frac{r_1 \Pi_{j\neq 1} s_j + \ldots + r_n \Pi_{j\neq n} s_j}{\Pi_1^n s_j};$$

$$\frac{\Pi_1^n s_j}{r_1 \Pi_{j\neq 1} s_j + \ldots + r_n \Pi_{j\neq n} s_j} (\frac{r_1}{s_1} p_1 + \ldots + \frac{r_n}{s_n} p_n)]$$

$$= [r_1/s_1 + \ldots + r_n/s_n; \frac{1}{r_1/s_1 + \ldots + r_n/s_n} \Sigma \kappa_i p_i] = [\Sigma_1^n \kappa_i; (1/\Sigma_1^n \kappa_i) \Sigma_1^n \kappa_i p_i].$$

In multi-valued PD, let $d_i = (p_{i1}, \ldots, _{im})$ for $i = 1, \ldots, n$ and $d = (p_1, \ldots, p_m)$ be distributions on a partition of $\Sigma$. Then $\alpha_i = [\kappa_i, d_i]$ confers credence $\kappa_i$ to all elements $p_{ij}$ of $d_i$. and $\alpha = [\kappa_1;d_1] \oplus, ..., \oplus [\kappa_n;d_n]$ is formed by merging $[\kappa_1;p_{1j}], ..., [\kappa_n;p_{nj}]$ on $U_2$ into $[\kappa;p_j]$ for $j=1, ..., m$, in accordance with rules (1) of binary mergers. Hence $\alpha = \{[\kappa_i;d_i]_1^n\} = [\kappa;d]$ with

$$\kappa = \Sigma_1^n \kappa_i; \quad p_j = (1/\Sigma_1^n \kappa_i) \Sigma_{i=1}^n \kappa_i p_{ij}; \quad j=1,...,m \tag{3.4}$$

An evidence with negative credence is 'counter-evidence' which cancels out with its positive counterpart, as in $[\kappa;p_1] \oplus [-\kappa;p_1] \oplus [\kappa_2;p_2] = [\kappa_2;p_2]$. Straight mergers are associative and commutative.

### 3.3.2 *Homomorphism*

A transformation is homomorphic if it transforms all functions of elements of a distribution into the same functions of the transformed elements.. Otherwise they are heteromorphic. Mergers are homomorphic if and only if they are associative and commutative. Heteromorphic mergers require rules that go beyond the rules of the probability calculus. Such rules are always somewhat arbitrary.

Straight mergers are homomorphic by equs. (4), since $\Sigma_1^n \kappa_i$ and $\Sigma_1^n \kappa_i p_{ij}$ are associative and commutative.

### 3.3.3 *Equivalence*

'$[\kappa;d]$ is true with probability $p$' is equivalent to $[p\kappa;d]$. (3.5)

*Proof*

Let $p=m/n$. Since $[\kappa;d] = \{[\kappa/n;d]|_1^n\}$, '$[\kappa;d]$ of probability $m/n$' is equivalent to $\{[\kappa/n;d]|_1^m\}$, which is $[m\kappa/n;d] = [p\kappa;d]$.

### 3.3.4 *Indifference*

PD expresses indifference to probabilities $p_1,...,p_n$ of a proposition by equi-credibility

$\kappa_1 = ... = \kappa_n$ of $[\kappa_1;p_1], ..., [\kappa_n;p_n]$.

## 3.4 Normalization

Let $U_m = \{A_1, ..., A_m\}$ be a partition of $\Sigma$, let $[\kappa_j;q_j]$ on $\{A_j, \overline{A}_j\}$ be a binary $\alpha$–evidence with $\kappa_j > 0$ and $q_j = P(A_j)$ the probability of $A_j$, and let us define $\beta = \{[\kappa_1;q_1], ..., [\kappa_m;q_m]\}$. PD interprets $\beta$ the 'normalized' commutative and associative $\alpha$–evidence $N\{\beta\} = [\hat{\kappa}; \hat{d}_m]$ with credence $\hat{\kappa}$ and distribution $\hat{d}_m = (\hat{p}_1, ..., \hat{p}_m)$ on $U_m$. $N\{\beta\}$ remains to be determined.

To determine $N\{\beta\}$ we note first that necessarily, if $(q_1,...,q_m) = (p_1,...,p_m) = d_m$ is a distribution and $\kappa_1 = ... = \kappa_m = \kappa$, then $N\{\beta\} = [\kappa;d_m]$.

$\alpha = [\kappa; d_m]$ on $U_m$ with $d_m = (p_1, \ldots, p_m)$ is equivalent to the set of binary evidences $\{[\kappa; p_1], \ldots, [\kappa; p_m]\}$ with $p_1 = P(A_1)$ etc., and by the rules of mergers,

$$\alpha = \{[\kappa p_j; \delta_{j1}, \ldots, \delta_{jj}, \ldots, \delta_{jm}]_1^m\}, \text{ where } \delta_{ij} = \begin{cases} 1 & \text{if } i = j; \\ 0 & \text{otherwise.} \end{cases}$$

This suggests an interpretation of $\beta$ as $\alpha$ – evidence

$$N\{\beta\} = \{[\kappa_j q_j; \delta_{j1}, \ldots, \delta_{jj}, \ldots, \delta_{jm}]_{j=1}^m\},$$

which satisfies the noted necessary condition. $N\{\beta\}$ is seen to be associative and commutative in $[\kappa_j; q_j]$. By Rules (3.1),

$$N\{\beta\} = [\Sigma \kappa_j q_j; \kappa_1 q_1 / \Sigma \kappa_j q_j, \ldots, \kappa_m q_m / \Sigma \kappa_j q_j; U_m],$$

i.e.,

*Rules of Normalization*

Norm $N\{\beta\} = [\hat{\kappa}; \hat{d}_m]$ with distribution $\hat{d}_m = (\hat{p}_1, \ldots, \hat{p}_m)$ on partition $U_m = \{A_1, \ldots, A_m\}$ of some sample space, of $\beta = \{[\kappa_1; q_1], \ldots, [\kappa_m; q_m]\}$ with $[\kappa_j; q_j]$ on $\{A_j, \overline{A}_j\}$ and $\kappa_j \geq 0$, is given by

$$\hat{\kappa} \hat{p}_j = \kappa_j q_j; \qquad \hat{\kappa} = \Sigma_1^m \kappa_j q_j. \qquad (3.6)$$

Probability $\hat{p}_j$ is the credence-weighted normalized value of $q_j$ and credence $\hat{\kappa}$ is a normalizing factor which enforces $\Sigma \hat{p}_j = 1$.

### 3.5    Offsetting Probability Dynamics (OPD)

In a stationary stochastic coin-tossing process, if half the outcomes in independent and equally credible successive trials are heads and the other half are tails, the result is taken as confirmation for probability ½ of heads, with increasing confidence the greater the number of trials. But if we accept equally credible reports concerning the outcome of a single toss and one half assert heads and the other half tails, then they should cancel out. Similarly, if in the epistemic theory two or more equi-credible and evidentially independent evidences for the same proposition have the same probability they reinforce each other, but if some assert 'true' while an equal number assert 'false', we have learnt nothing and they should cancel out. Offsetting probability dynamics (OPD) is a generalization of these examples.

$\alpha_1 \diamond \ldots \diamond \alpha_n$ is the OPD-merger of evidentially independent $\alpha_i, \ldots, \alpha_n$, where $\alpha_1 = [\kappa_1; p_1]$, etc. OPD merger $[\kappa; 0] \diamond [\kappa; 1]$ of the two binary equi-credible evidences $\alpha_1 = [\kappa; 0]$ and $\alpha_2 = [\kappa; 1]$ should have credence nil, and $n$ evidentially independent binary evidences $[\kappa_1; p], \ldots, [\kappa_n; p]$ asserting the same $p$ with respective credence $\kappa_1,\ldots,\kappa_n$ should merge in OPD as they do in SPD.

OPD merger $\breve{\alpha}(n)=[\breve{\kappa}(n);\breve{p}(n)]=$'$\alpha_1 \diamondsuit \ldots \diamondsuit \alpha_n$' is defined by the following postulates:-

*Postulates of OPD-Mergers*

OP1     In $\breve{\alpha}(n)=\alpha_1 \diamondsuit \ldots \diamondsuit \alpha_n = \breve{\kappa}(n);\breve{p}(n)$, where $\alpha_i=[\kappa_i;p_i]$; i$i$=1,…,$n$; increasingly conflicting distributions ($p_1,\ldots,p_n$) increasingly offset each other. In particular, if $n$ is even then $\alpha_i=[\kappa_i;1]$ for $i=1,\ldots,\tfrac{1}{2}n$, $\alpha_i=[\kappa_i;0]$ for $i=\tfrac{1}{2}n+1,\ldots,n$, and $\Sigma_1^{\tfrac{1}{2}n}\kappa_i=\Sigma_{\tfrac{1}{2}+1}^n\kappa_i$, then $\breve{\kappa}(n)=0$.

OP2     If $p_1=\ldots=p_n$ then OPD is SPD.

OP3     If $\kappa_i=0$ then $\breve{\alpha}(n)$ is independent of $\alpha_i$.

OP4     The elements of $\breve{\alpha}(n)$ are continuous functions of the elements of the merged $\alpha_i$.

OP5     Scaling the credence of all merged evidences by a common factor scales the credence of $\breve{\alpha}(n)$ by the same factor and does not affect its probability.

OP6     $[\breve{\kappa}(n);\breve{p}(n)]=[\kappa_1;p_1]\diamondsuit \ldots \diamondsuit [\kappa_n;p_n]$ is equivalent to
$[\breve{\kappa}(n);1-\breve{p}(n)]=[\kappa_1;1-p_1]\diamondsuit \ldots \diamondsuit [\kappa_n;1-p_n]$.

OP7     $\breve{\alpha}(n)$ is associative and commutative in the merged $[\kappa_i;p_i]$.

Assuming $\kappa_i \geq 0$, let us determine a measure $\lambda(n)$ of the 'accord' of $[\kappa_1;p_1],\ldots,[\kappa_n;p_n]$ on U$_2$, such that $\breve{\alpha}(n)=[\breve{\kappa}(n);\breve{p}(n)]$ with

$$\breve{p}(n)=\overline{p}(n)=\Sigma_1^n \kappa_i p_i / \Sigma_1^n \kappa_i; \quad \breve{\kappa}(n)=\lambda(n)\Sigma_1^n \kappa_i; \qquad (3.7)$$

satisfies the postulates of OPD mergers.

According to eqs. (7), $\breve{\alpha}(n)$ has the probability of a straight merger but credence $\lambda(n)\Sigma_1^n \kappa_i$ rather than $\Sigma_1^n \kappa_i$. Accord $\lambda(n)$ is total and has $\lambda(n)=1$ if and only if $p_1=\ldots=p_n$. 'Discord' $1-\lambda(n)$ is total if and only if $\lambda(n)=0$. $\lambda(n)$ decreases with an appropriately defined credence-weighted dispersion of elements $[\kappa_1;p_1],\ldots,[\kappa_n;p_n]$, from 1 when there is no dispersion to 0 when the evidences jointly annul each other.

$\breve{\alpha}(n)=[\breve{\kappa}(n);\breve{p}(n)]=[\kappa_1;p_1]\diamondsuit \ldots \diamondsuit [\kappa_n;p_n]$, with

$$\breve{p}(n)=\overline{p}(n)=\Sigma_1^n \kappa_i p_i / \Sigma_1^n \kappa_i; \quad \breve{\kappa}(n)=\lambda(n)\Sigma_1^n \kappa_i; \qquad (3.8)$$

yields OPD-mergers that satisfy the Postulates if $\lambda(n)$ is given by

$$\lambda(n) = 1 - 2\sigma(n);$$

$$[\sigma(n)]^2 = \Sigma_1^n \kappa_i [p_i - \overline{p}(n)]^2 / \Sigma_1^n \kappa_i; \quad \sigma(n) \geq 0; \quad \quad (3.9)$$

$$\overline{p}(n) = \Sigma_1^n \kappa_i p_i / \Sigma_1^n \kappa_i.$$

$[\sigma(n)]^2$ and $\lambda(n)$ remain invariant if every $p_i$ in these equations is replaced by $1 - p_i$. $\sigma(n)$ is a measure of the dispersion of the random variable $p$ that ranges over $(p_1,...,p_n)$ with distribution $(\kappa_1 / \Sigma \kappa_i, ..., \kappa_n / \Sigma \kappa_i)$. $\sigma(n)$ is the standard deviation of $p$ and $\overline{p}(n) = Exp[p]$ is its average. Clearly, $\lambda(n) \leq 1$. It can be shown that $\lambda(n) \geq 0$.

In case of equi-credibility, $\kappa_1 = ... = \kappa_n = \kappa$,

$$\lambda(n) = 1 - 2\sqrt{\tfrac{1}{n}[\Sigma p_i^2 - (\Sigma p_i)^2]},$$

which is independent of credence, and

$$\breve{p}(n) = \overline{p}(n) = \tfrac{1}{n}\Sigma p_i; \quad \breve{\kappa}(n) = \lambda(n)n\kappa.$$

If $n=2$, equs. (9) reduce to $\lambda(2) = \lambda$ with

$$\lambda = 1 - 2|p_1 - p_2|\sqrt{\kappa_1 \kappa_2}/(\kappa_1 + \kappa_2). \quad \quad (3.10)$$

### 3.6 Correlation Coefficients

Correlation coefficient $\rho(A, B)$ is a measure of the mutual interaction of A and B. It is determined by $E = [P(A), P(AB), P(B)]$ and defined by

$$\rho(A,B) = Exp[(A - \mu)(B - \nu)]/\sigma_A \sigma_B = [Exp\, AB - \mu\nu]/\sigma_A \sigma_B.$$

Here $Exp$ = mathematical expectation, $\mu = Exp[A]$, $\nu = Exp[B]$, and

$$\sigma_A = \sqrt{Exp[(A - \mu)^2]} = \sqrt{Exp A^2 - \mu^2}; \quad \sigma_B = \sqrt{Exp B^2 - \nu^2};$$

are standard deviations. If A and B are propositions, they are binary, true or false. In this case

$$\mu = P(A), \quad \nu = P(B), \quad Exp[AB] = P(AB),$$

$$\sigma_A^2 = P(A)P(\overline{A}); \quad \sigma_B^2 = P(B)P(\overline{B});$$

$$\rho(A, B) = [P(AB) - P(A)P(B)] / \sqrt{P(A)P(\overline{A})P(B)P(\overline{B})}. \quad \quad (3.11)$$

$\rho(A, B)$ corresponds to $E$, $\rho_0(A, B)$ corresponds to $E_0$, $\rho_{10}(A, B)$ corresponds to $E_{10}$, etc.

$$-1 \leq \rho(A, B) \leq 1.$$

$$\rho(A, B) = \rho(B, A) = -\rho(A, \overline{B}) = -\rho(\overline{A}, B).$$

$\rho(A, B) \to 0$ if and only $P(A|B) \to P(A)$. [8]

If $P(A)$ or $P(B) = 0$ or $1$ then $\rho(A, B) = 0$.

$$\lim_{P(A) \to 0} \rho(A, B) = \lim_{P(A) \to 1} \rho(A, B) = \lim_{P(B) \to 0} \rho(A, B) = \lim_{P(B) \to 1} \rho(A, B) = 0. \qquad (3.12)$$

We prove that $\lim_{P(B) \to 0} \rho(A, B) = 0$ if $0 < P(A) < 1$:-

$$\lim_{P(B) \to 0} \rho(A, B) = \lim_{P(B) \to 0} [P(AB) - P(A)P(B)] / \sqrt{P(A)P(\overline{A})P(B)P(\overline{B})}$$

$$= \lim_{P(B) \to 0} [\sqrt{P(B)}P(AB))/P(B) - P(A)\sqrt{P(B)}] / \sqrt{P(A)P(\overline{A})} = 0.$$

All other relations of equ.(12) now follow at once.

### 3.7 Cross-Credence[9]

For non-negative $\kappa_1$ and $\kappa_2$ let $[\kappa; d]$ result from accepting $[\kappa_1; d]$ with credence $\kappa_2$. Then 'cross-credence' $\kappa$ of $\kappa_1$ and $\kappa_2$ should satisfy the following conditions:

C1     $\kappa = f(\kappa_1, \kappa_2)$ for some continuously differentiable commutative function $f(.)$ of $\kappa_1$ and $\kappa_2$.
C2     If $\kappa_1$ and $\kappa_2$ are scaled by a common factor then $\kappa$ is scaled by the same factor.
C3     $f(0,x) = f(x,0) = 0$ for all real x.
C4     $f(\infty, x) = x$ for all real x.
C5     $f(x,x) = \frac{1}{2}x$.

Condition C1 is evidently reasonable and C2-C4 are clearly necessary but C5 is only justified by considerations of simplicity. C1-C4 entail '$\kappa^\nu = \kappa_1^\nu + \kappa_2^\nu$' for some negative $\nu$. We set $\nu = -1$ because we regard it the simplest choice. Taken together with C1-C4, C5 is equivalent to this choice. Thus

$$1/\kappa = 1/\kappa_1 + 1/\kappa_2, \qquad (3.13)$$

$$\kappa = \kappa_1 \kappa_2 / (\kappa_1 + \kappa_2). \qquad (3.14)$$

With non-negative $\kappa_1,...,\kappa_n$, let the proposition '''$[\kappa_1;d]$ is of credence $\kappa_2$' is of credence $\kappa_3$' etc. up to $\kappa_n$', i.e., '$[\kappa_n;[...[\kappa_2;[\kappa_1;d]]...]]$', be equivalent to $[\kappa;d]$. Then $\kappa$ is the cross-credence of $\kappa_1,...,\kappa_n$, and it follows from equ. (13) that

*Theorems of Cross-Credence*

Cross-credence $\kappa$ of non-negative $\kappa_1,...,\kappa_n$ is given by

$$1/\kappa = \Sigma_{i=1}^{n} 1/\kappa_i; \qquad \kappa_i \geq 0. \tag{3.15}$$

Next, let $d_1$ and $d_2$ be distributions of probabilities of A and B, respectively, with respective credences $\kappa_1$ and $\kappa_2$, let $[\kappa_1;d_1]$ be given, and let the acceptance of evidence $[\kappa_2;d_2]$ with credence $\kappa_1$ produce $[\kappa;d_2]$. Then $\kappa$ is the unilateral cross-credence from A to B of $\kappa_1$ and $\kappa_2$ of equ. (14), except that just that fraction of $\kappa_1$ is effective in the formation of $\kappa$ which expresses the degree of interaction of A and B. An adequate expression for this fraction is the magnitude of correlation coefficient $\rho(A,B)$ of A and B. This transforms $\kappa$ of equ. (14) into

$$\kappa^* = \frac{|\rho(A,B)|\kappa_1\kappa_2}{|\rho(A,B)|\kappa_1 + \kappa_2}. \tag{3.16}$$

Substitution in equ. (14) of $|\rho(A,B)|\kappa_1$ for $\kappa_1$ and $|\rho(A,B)|\kappa_2$ for $\kappa_2$ yields the bilateral cross-credence

$$\kappa^{**} = |\rho(A,B)|\frac{\kappa_1\kappa_2}{\kappa_1 + \kappa_2}. \tag{3.17}$$

**3.8　Updating Fallible Evidence**

On partition $\{A,\overline{A}\} \times \{B_1,...,B_m\}$, given prior

$$E_0 = \{P_0(AB_1),...,P_0(AB_m)\} \tag{3.18}$$

at $t_0$, let us determine update $P_{10}(A)$ of $P_0(A)$ by distribution $\tilde{d}_1^B = (\tilde{P}_1(B_1),...,\tilde{P}_1(B_m))$ accepted at $t_1 > t_0$. Writing $P_{10}(A)$ for $P_1(A)$ of equ. (1.5), then, if $\tilde{P}_1(B_1) = 1$, acceptance of $\tilde{d}_1^B$ produces $P_{10}(A) = P_0(A|B)$. Therefore acceptance of $\tilde{d}_1^B$ updates $P_0(A)$ into $P_0(A|B_1)$ with probability $\tilde{P}_1(B_1)$, into $P_0(A|B_2)$ with probability $\tilde{P}_1(B_2)$, etc., i.e., acceptance of $\tilde{d}_1^B$ updates $P_0(A)$ into the mathematical expectation of the conditional probability $P_0(A|B)$ of equ. (1.5) with B ranging over $\{B_1,...,B_m\}$ with distribution $\tilde{d}_1^B$. Thus

$$P_{10}(A) = \Sigma_1^m \tilde{P}_1(B_j) P_0(A|B_j), \tag{3.19}$$

whose binary instance is

$$P_{10}(A)=\widetilde{P}_1(B)P_0(A|B)+\widetilde{P}_1(\overline{B})P_0(A|\overline{B}). \tag{3.20}$$

These are the defining rules of 'Probability Kinematics' (Jeffrey, 1965) which go some way toward a dynamic theory. Yet Probability Kinematics is inconsistent. As an infallible $P_0(A)$ cannot be updated, it tacitly assumes that $P_0(A)$ is fallible, but it disregards its degree of fallibility, and it is incomplete because it accepts $\widetilde{P}_1(B_j)$ at $t_1$ with certainty and thereby equates it with update $P_1(B_j)$ of $P_0(B_j)$, whereas a higher order theory should be able to admit fallible $\widetilde{P}_1(B_j)$ and distinguish them from $P_1(B_j)$.

Update $P_{10}(A)$ of $P_0(A)$ produced by Probability Kinematics is the weighted sum of conditional first order probabilities $P_0(A|B_1),\ldots,P_0(A|B_m)$, with $\widetilde{P}_1(B_1),\ldots,\widetilde{P}_1(B_m)$ as respective weights, in contradistinction to the second order probabilities discussed in Section 2.2, which are probabilities of probabilities.

In Probability Dynamics, prior[10] $\alpha_0^{AB}=[\kappa_0;d_0^{AB}]$ has distribution

$$d_0^{AB}=(P_0(AB_1),\ldots,P_0(AB_m),P_0(\overline{A}B_1),\ldots,P_0(\overline{A}B_m)), \tag{3.21}$$

whose binary instance is

$$d_0^{AB}=((P_0(AB),P_0(A\overline{B}),P_0(\overline{A}B),P_0(\overline{A}\,\overline{B})). \tag{3.22}$$

If $\alpha_0^{AB}$ is given, then so are $\alpha_0^A=[\kappa_0;P_0(A),P_0(\overline{A})], \alpha_0^{B_j}=[\kappa_0;P_0(B_j),P_0(\overline{B}_j)]$ and correlation coefficients $\rho_0(A,B_j)$.

Let us define $\alpha_1^{AB}=[\kappa_1,d_1^{AB}]$ as $\alpha_0^{AB}$ with $\kappa_1$ for $\kappa_0$ and $d_1^{AB}$ for $d_0^{AB}$ with $P_1$ for $P_0$. Then prior $\alpha_0^{AB}$ is updated by $\alpha_1^{AB}$ into their merger.

Next, we determine on $\{A,\overline{A}\}\times\{B_1,\ldots,B_m\}$ the 'indirect' update

$$\alpha_{10}=[\kappa_{10};d_{10}]; \quad d_{10}=(P_{10}(A),P_{10}(\overline{A})), \tag{3.23}$$

of

$$\alpha_0=[\kappa_0;d_0]; \quad d_0=(P_0(A_1),P_0(\overline{A})), \tag{3.24}$$

by evidence

$$\widetilde{\alpha}_1^B=[\widetilde{\kappa}_1;\widetilde{d}_1]; \quad \widetilde{d}_1=(\widetilde{P}_1(B_1),\ldots,\widetilde{P}_1(B_m)) \tag{3.25}$$

and nothing else accepted at $t_1$, given prior $\alpha_0^{AB}$ at $t_0$ with $d_0^{AB}$ of equ. (21).

For this purpose we first form, for $j=1,\ldots,m$, the binary update $\alpha_{10}^j=[\kappa_{10}^j;P_{10}^j(A),P_{10}^j(\overline{A})]$ of $[\kappa_0;P_0(A),P_0(\overline{A})]$ by $[\widetilde{\kappa}_1;\widetilde{P}_1(B_j),\widetilde{P}_1(\overline{B}_j)]$. Then $\kappa_{10}^j$ is the unilateral cross-credence $\kappa^*$, equ. (3.16), with $\widetilde{\kappa}_1$ for $\kappa_1$, $\kappa_0$ for $\kappa_2$, and $B_j$ for B, i.e.,

$$\kappa_{10}^{j} = \frac{|\rho(A,B_j)|\,\widetilde{\kappa}_1 \kappa_0}{|\rho(A,B_j)|\,\widetilde{\kappa}_1 + \kappa_0}. \qquad (3.26)$$

In general, $\kappa_{10}^j$ depends on $j$ and set $\{[\kappa_{10}^1;\widetilde{P}_1(B_1)],\ldots,[\kappa_{10}^m;\widetilde{P}_1(B_m)]\}$ is no $\alpha$-evidence. Therefore we normalize it into $\alpha$-evidence $\hat{\alpha}_{10}=[\hat{\kappa}_{10};\hat{P}_1(B_1),\ldots,\hat{P}_1(B_m)]$, and obtain by Rules (6)

$$\hat{\kappa}_{10} = \Sigma_1^m \kappa_{10}^j \widetilde{P}_1(B_j); \quad \hat{P}_1(B_j) = (1/\hat{\kappa}_{10})\,\kappa_{10}^j\,\widetilde{P}_1(B_j). \qquad (3.27)$$

Update $\alpha_{10}$, equ. (23), of $\alpha_0^A$, equ. (24), by $\widetilde{\alpha}_1^B$, equ. (25), is now produced with credence $\kappa_{10} = \hat{\kappa}_{10}$ and distribution $d_{10} = (P_{10}(A), P_{10}(\overline{A}))$ of equ. (23) and $P_{10}(A)$ is the mathematical expectation $Exp[P_0(A|B)]$ of $P_0(A|B)$ with B ranging over $\{B_1,\ldots,B_m\}$ with distribution $(\hat{P}_1(B_1),\ldots,\hat{P}_1(B_m))$, i.e.,

$$\alpha_{10} = [\hat{\kappa}_{10}; P_{10}(A), P_{10}(\overline{A})]; \quad P_{10}(A) = \Sigma_1^m \hat{P}_1(B_j)\,P_0(A|B_j). \qquad (3.28)$$

In the binary case ($m=2$) we set $B_1 = B$, $B_2 = \overline{B}$, and have $|\rho_0(A,B)| = |\rho_0(A,\overline{B})|$ and $\kappa_{10}^1 = \kappa_{10}^2 = \kappa_{10}$. Therefore this case requires no normalization. Update $\alpha_{10}$, equ. (23), of $\alpha_0$, equ. (24), by $\widetilde{\alpha}_1^B = [\widetilde{\kappa}_1; \widetilde{P}_1(B), \widetilde{P}_1(\overline{B})]$ is now given by

$$\kappa_{10} = \frac{|\rho_0(A,B)|\,\widetilde{\kappa}_1 \kappa_0}{|\rho_0(A,B)|\,\widetilde{\kappa}_1 + \kappa_0}; \qquad (3.29)$$

and

$$P_{10}(A) = \widetilde{P}_1(B) P_0(A|B) + \widetilde{P}_1(\overline{B}) P_0(A|\overline{B}). \qquad (3.30)$$

This solves the problem of the production of update $\alpha_{10}$ of $\alpha_0$ by $\widetilde{\alpha}_1^B$.

In SPD, prior $\alpha_0$ is successively updated by $\widetilde{\alpha}_1^B = [\widetilde{\kappa}_1; \widetilde{d}_1^B], \ldots, \widetilde{\alpha}_n^B = [\widetilde{\kappa}_n; \widetilde{d}_n^B]$ by first forming merger $[\widetilde{\kappa};\widetilde{d}] = \widetilde{\alpha}_1^B \oplus \ldots \oplus \widetilde{\alpha}_n^B$ and subsequently updating $\alpha_0$ by $[\widetilde{\kappa};\widetilde{d}]$. Since $[\widetilde{\kappa};\widetilde{d}]$ is commutative in $\widetilde{\alpha}_1^B,\ldots,\widetilde{\alpha}_n^B$, so is the resulting update, in contrast to the updating process of Probability Kinematics which is not commutative but is claimed to have commutative variants [cf. Field, 1978; Lange, 2000; Wagner, 2002, 2003; Hawthorne, 2004].

Probability Kinematics emerges from Probability Dynamics as $\widetilde{\kappa}_1 \to \infty$, but this limiting case makes update $P_{10}(B)$ of $P_0(B)$ identical to $\widetilde{P}_1(B)$, whereas the two should generally be distinct. It follows that Probability Kinematics is unsound.

If both $\widetilde{\kappa}_1$ and $\widetilde{\kappa}_2$ tend to infinity then $[\widetilde{\kappa};\widetilde{p}] = [\widetilde{\kappa}_1;\widetilde{p}_1] \oplus [\widetilde{\kappa}_2;\widetilde{p}_2]$ is indefinite and depends on how this is brought about. For example, if $n=2$ and $\widetilde{\kappa}_1 \to \widetilde{\kappa}_2$ as $\tilde{\kappa}_1 \to \infty$ then $\tilde{p} \to \frac{1}{2}(\tilde{p}_1 + \tilde{p}_2)$, but if $\tilde{\kappa}_1 \to \frac{2}{3} \tilde{\kappa}_2$ as $\tilde{\kappa}_1 \to \infty$ then $\tilde{p} \to \frac{2}{3} \tilde{p}_1 + \frac{1}{3} \tilde{p}_2$. Therefore PD does not confer its commutativity to Probability Kinematics.

The same conclusions follow from the observation that the acceptance of $[\tilde{\kappa}; \tilde{P}(B)]$ with credence $\tilde{\kappa} = \infty$ means an acceptance of $\tilde{P}(B)$ with certainty, but both $\tilde{P}_1(B)$ and $\tilde{P}_2(B)$ cannot be certain unless they are equal.

### 3.9   Probability Dynamics and Truth

PD is a theory whose object is the provision of estimates of probabilities from fallible evidence. Epistemic PD need not be more than coherent and reasonable. The verity of epistemic PD is relative to its assumptions. For example, criteria of optimality are always somewhat arbitrary.

Frequentist PD is applied in descriptions of reality. It sets out with probabilities and credence resulting from observation reports, and its conclusions should be as free of arbitrary assumptions as possible.

### PART IV         New Explanation of Old Evidence [11, 12]

### 4.1   *The Problem*

On $\Omega = \{A, \overline{A}\} \times \{B, \overline{B}\}$ with A an 'hypothesis' and B an 'evidence' (not to be confounded with evidences $E$ and $\alpha$ of Parts I and III), let

$$\alpha = \begin{pmatrix} a & b \\ c & d \end{pmatrix} \begin{matrix} \text{B} & \overline{\text{B}} \\ \text{A} \\ \overline{\text{A}} \end{matrix} \qquad (4.1)$$

be a distribution with probabilities $a = P(AB)$, $b = P(A\overline{B})$, $c = P(\overline{A}B)$, $d = P(\overline{A}\,\overline{B})$.

Jeffrey (1995) discusses the problem of the transformation of such an initially given 'defective' distribution $\alpha$ into a 'repaired' distribution

$$\breve{\alpha} = \begin{pmatrix} \breve{a} & \breve{b} \\ \breve{c} & \breve{d} \end{pmatrix}, \qquad (4.2)$$

by the subsequent independent discovery that A implies B. In this 'problem of old evidence and new explanation' he obtains by methods of first order probability and probability kinematics, that

$$\breve{\alpha} = \begin{pmatrix} a+b & 0 \\ c & d \end{pmatrix}. \qquad (4.3)$$

The more general case of probable old evidence and new explanation is the subject of several papers by Wagner (1997 et seq.) which are not on the right track.

If $b=0$ in $\alpha$ of equ. (1) then A implies B, and if $b>0$ then $\alpha$ is not compatible with an implication of B by A, and it is possible to proceed only if $\alpha$ is fallible, but this is not within the scope of a first order theory.

We apply the methods of PD and obtain results which do not sustain those of Jeffrey and his followers.

### 4.2 *The Solution*

Let 'A→B' be a conditional implication which implies B if A is true and is evidentially neutral toward B otherwise. Let evidence

$$\beta = [\kappa;\alpha] = [\kappa;\begin{pmatrix} a & b \\ c & d \end{pmatrix}] \qquad (4.4)$$

be given and let us accept 'constraint'

$$\widetilde{\beta} = [\widetilde{\kappa};P(A\rightarrow B)]. \qquad (4.5)$$

Then, given $\beta$, an acceptance of $\widetilde{\beta}$ produces

$$\breve{\beta} = [\breve{\kappa};\breve{\alpha}] = [\breve{\kappa};\begin{pmatrix} a+b & 0 \\ c & d \end{pmatrix}]; \quad \breve{\kappa} = \widetilde{\kappa}\kappa/(\widetilde{\kappa}+\kappa). \qquad (4.6)$$

*Proof*

The imposition of 'A→B' on $\alpha$ 1) does not affect $c$ and $d$ since, by definition, 'A→B' is inactivated if B is false, 2) transforms $b$ into 0, because 'A→B' is true if A is true if and only if B is true, and 3) replaces $a$ by $a+b$ since $(a+b)+c+d=1$. Therefore $\breve{\alpha}$ is given by equ. (3). $\breve{\kappa}$ is cross-credence $\widetilde{\kappa}\kappa/(\widetilde{\kappa}+\kappa)$ of $\widetilde{\kappa}$ and $\kappa$, which completes the proof.

The acceptance of $\breve{\beta}$ revises $\beta$ into their merger

$$\beta_1=[\kappa_1;\alpha_1]=\beta\oplus\breve{\beta}=[(\kappa+\breve{\kappa});(1/(\kappa+\breve{\kappa}))(\kappa\alpha+\breve{\kappa}\breve{\alpha})]. \qquad (4.7)$$

This solves our problem.

### 4.3 *Discussion*

*a)* $\beta_1$ of equ. (7) retains $\alpha$ (as part of $\beta$) in contrast with first order theories which cannot simultaneously accommodate both $\alpha$ and $\widetilde{\alpha}$ and are thereby driven to abandon $\alpha$, although they continue to make use of its elements. The root of the matter is the inability of first order theories to form mergers.

*b)* If $(\widetilde{\kappa}/\kappa)\rightarrow\infty$ then $\breve{\kappa}\rightarrow\kappa$, $\kappa_1\rightarrow 2\kappa$, and $\alpha_1\rightarrow(\alpha+\breve{\alpha})$. $(\widetilde{\kappa}/\kappa)\rightarrow 0$

inactivates the imposition of constraint $\widetilde{\beta}$ of equ. (5) and causes $\beta_1$ to tend toward $\beta$, as it should, because ($\widetilde{\kappa}/\kappa$)=0 means that no constraints are imposed on $\alpha$. The imposition of [$\widetilde{\kappa}$;P(A→B)] on $\alpha$ retains $\alpha$ when it produces $\beta_1$: The problem is not how to repair a defective distribution but how to incorporate new evidence.

c)  In equ. (7) $\alpha_1$ is

$$\alpha_1 = (1/\kappa_1)\begin{pmatrix} (\kappa+\widetilde{\kappa})a+\widetilde{\kappa}b & \kappa b \\ (\kappa+\widetilde{\kappa})c & (\kappa+\widetilde{\kappa})d \end{pmatrix}; \quad \kappa_1 = \kappa + \breve{\kappa} . \tag{4.8}$$

The resulting probability of A conditional on B is

$$P_1(A|B) = P_1(AB)/P_1(B) = \frac{a+(\widetilde{\kappa}/(\kappa+\widetilde{\kappa}))b}{a+(\widetilde{\kappa}/(\kappa+\widetilde{\kappa}))b+c} . \tag{4.9}$$

$P_1(A|B)$ increases with $\widetilde{\kappa}/\kappa$ from $P_1(A|B)=P(A|B)=a/(a+c)$ at ($\widetilde{\kappa}/\kappa$)=0 up to $(a+b)/(a+b+c)$ at ($\breve{\kappa}/\kappa$)=∞. If ($\breve{\kappa}/\kappa$)>0 then $P_1(A|B)$>$P(A|B)$. Acceptance of constraint $\widetilde{\beta}$ is seen to increase P(A|B).

d)  According to the theorem of equivalence, equ. (3.5), acceptance of constraint [$\widetilde{\kappa}$;P(A→B)=$p$] instead of $\widetilde{\beta}$ of equ. (5), merely replaces $\widetilde{\kappa}$ in [$\widetilde{\kappa}$;P(A→B)] by $p\widetilde{\kappa}$.

## 4.4  *A Generalization*

Finally, let us define on $\Omega=\{A,\overline{A}\}\times\{B_1,...,B_n\}$ evidences $\beta=[\kappa;\alpha]$ and $\beta_i=[\kappa_i;\alpha_i]$ for $i=1,...,n$, with distributions

$$\alpha = \begin{pmatrix} a_{11}...a_{1i}...a_{1n} \\ a_{21}...a_{2i}...a_{2n} \end{pmatrix}; \quad \alpha_i = \begin{pmatrix} 0...0 \ \Sigma a_{1j} \ 0...0 \\ a_{21} \ ... \ a_{2i} \ ... \ a_{2n} \end{pmatrix}. \tag{4.10}$$

To keep the exposition simple we restrict the sequel to straight mergers. It then follows by the same reasoning as before and a little algebra that the successive imposition on $\beta$ of $n$ constraints $C_i$=[$\kappa_i$;P(A→$B_i$)] with $i$=1,…,$n$, generates the commutative evidence

$$\breve{\beta} = [\breve{\kappa}_1;\alpha_1] \oplus ... \oplus [\breve{\kappa}_n;\alpha_n] = [\Sigma\breve{\kappa}_i;\begin{pmatrix} \mu\kappa_1 \ ... \ \mu\kappa_n \\ b_1 \ ... \ b_n \end{pmatrix}]; \tag{4.11}$$

with $\mu=\Sigma a_i/\Sigma\breve{\kappa}_i$ by definition and $\breve{\kappa}_i=\widetilde{\kappa}_i\kappa/(\widetilde{\kappa}_i+\kappa)$; $i=1,...,n$. An imposition of constraints $C_i$ with respective probabilities $p_i$, replaces every $\kappa_i$ by $p_i\kappa_i$. The resulting revision of $\beta$ is $\beta_1=\beta\oplus\breve{\beta}$.

**PART V   PDCT:  The Theory of Confirmation**

The first order theory defines the 'inductive support' or 'degree of confirmation' of 'hypothesis' A by 'evidence' B accepted at $t_1$, as excess $C_{10}$ of the updated probability $P_1(A)$ of A over its prior probability $P_0(A)$, i.e.,

$$C_{10} = P_1(A) - P_0(A), \tag{5.1}$$

and by equ. (1.5), given $E_0$ of equ. (1.4),

$$C_{10} = P_0(A|B) - P_0(A). \tag{5.2}$$

PD provides two different definitions of degrees of confirmation. The first defines the degrees of confirmation '$(\kappa_{10} - \kappa_0)/\kappa_0$' for credence and '$P_{10}(A) - P_0(A)$' for probability. The second defines degree of confirmation $K_{10}(A,B)$ of $\alpha_0 = [\kappa_0; P_0(A)]$ by $\tilde{\alpha}_1^B = [\tilde{\kappa}_1; \tilde{P}_1(B)]$ given prior $\alpha_0^{AB} = [\kappa_0; d_0^{AB}]$, with $d_0^{AB}$ of equ. (3.22), as the excess of the product of credence and probability of update $\alpha_{10} = [\kappa_{10}; P_{10}(A)]$ over that of prior $\alpha_0$, both divided by $\kappa_0$ as reference credence.

In this second definition, if $\alpha_0$ is updated by accepting $\tilde{\alpha}_1^B$ and nothing else, then

$$K_{10}(A,B) = \lambda \, [\kappa_{10}/\kappa_0] P_{10}(A) - (\kappa_0/\kappa_0) P_0(A),$$

with $\lambda$ the measure of accord of Section 3.5, and by (3.29) and (3.30),

$$K_{10}(A,B) = \lambda \, \frac{|\rho_0(A,B)| \, \tilde{\kappa}_1}{|\rho_0(A,B)| \, \tilde{\kappa}_1 + \kappa_0} [\tilde{P}_1(B) P_0(A|B) + \tilde{P}_1(\overline{B}) P_0(A|\overline{B})] - P_0(A). \tag{5.3}$$

In straight PDCT $\lambda = 1$, and if $(\tilde{\kappa}_1/\kappa_0) = \infty$ and $\tilde{P}_1(B) = 1$, then $K_{10}(A,B)$ reduces to $C_{10}$ of equ. (1). In binary offsetting PDCT, by equ. (3.10) with $P_0(A)$ for $p_1$, $P_{10}(A)$ for $p_2$, $\kappa_0$ for $\kappa_1$, and $\tilde{\kappa}_1$ for $\kappa_2$,

$$\lambda = 1 - 2|P_0(A) - P_{10}(A)| \sqrt{\kappa_0 \tilde{\kappa}_1} / (\kappa_0 + \tilde{\kappa}_1). \tag{5.4}$$

**Part VI – The Principle of Reflection**

**6.1** *The Principle*

The so-called 'Principle of Reflection' (Van Fraassen, 1984) has brought about many discussions pro and contra (e.g. Van Fraassen, 1989, 1999 [13]; Gaifman, 1986; Maher, 1992; Kvanvig, 1994; Green et al., 1994; Hild, 1998; Loewer, 1998, and many others). It consists of a verbal and of a mathematical principle. We show that the verbal principle does not follow from the mathematical version, prove that both are false, and analyse Van Fraassen's attempted proof.

The mathematical version purports to make second order theories of probability consistent by imposing the constraint discussed in Section 2.4. Our discussion illustrates pitfalls of inadequate formalizations and of mistaken conclusions from so-called 'Dutch book' – heads I win, tails you lose – betting schemes.

According to the *Verbal Principle of Reflection*

R1  "[An] agent's present subjective probability for proposition A, on the supposition that his subjective probability for this proposition will equal *r* at some later time, must be equal to this same number *r*."

Van Fraassen claims that this principle is equivalent to

R2  *The Special Mathematical Principle of Reflection,*

"$P_t^a[A \mid p_{t+x}^a(A) = r] = r$, where $P_t^a$ is agent a's credence function at time *t*, *x* is any non-negative number, and ($p_{t+x}^a(A) = r$) is the proposition that at time *t+x* the agent will bestow degree *r* of credence on the proposition A."

R2 is ambiguous, but if we read $P_t^a[A \mid p_{t+x}^a(A) = r]$ as an agent's probability of A at time *t* conditional on its value *r* at time $t + x$, and introduce a prior, we obtain for it in my notation

R3  Given $P_0[P_1(A) = r] = 1$ and prior $E_0$, then $P_0[A \mid P_1(A) = r] = r$,

i.e., at $t_0$, given $E_0$, an agent's subjective probability of A conditional on his subjective probability of A being *r* at some later time $t_1$, is equal to *r*. R3 is a special case of R4.

If $P_t^a[A \mid p_{t+x}^a(A) = r]$ were meant to say that $P_t^a$ is the probability of A at time *t given* '$p_{t+x}^a(A) = r$' (our '$P_1(A)=r$'), it would not be a conditional probability but probability $P_0[P_1(A)=r. \rightarrow A]$ of the conditional, which vanishes identically if *r*<1 since $[P_1(A)=r. \rightarrow A]$ is true if and only if *r*=1.

The probability of valid propositions (tautologies) is one, but probabilities of epistemic propositions are neither nil nor one. Therefore conditionalization is of little interest if $P_0[P_1(A) = r] = 1$.

R4    *The Mathematical Principle of Reflection*

Given $E_0 = $ '$P_0[P_1(A) = r] \& P_0[A \wedge .P_1(A) = r]$' and $0 < P_0[P_1(A) = r] < 1$, then

$$P_0[A | P_1(A) = r] = r.$$

Van Fraassen is right to devote to it the greater part of his considerations to R4, although this is not the way he puts it.

## 6.2    *Discussion*

R4 is trivially true if $r=0$ or 1. Therefore we exclude these instances from our discussion and assume that $0<r<1$. For the moment let us disregard that $P_0[P_1(A) = r]$ is always nil for epistemic A (Section 2.2).

In R2, A is conditional on ($p^a_{t+x}(A) = r$) and in R3 it is conditional on '$P_1(A) = r$'. In R1 the 'supposition' that the probability of A "will equal $r$ at some later time" is ambiguous. If '$p^a_{t+x}(A) = r$' is the conditional $P_0(A|B)$ then it is not the agent's (unconditional) probability of A as claimed in R1, and if it is the indicative $P_1(A | B)$ it is the agent's probability of A at $t_1$ and not at $t_0$. In both cases it has nothing to do with reflections. Nor does R1 follow from R2 and R3. A fortiori it does not follow from R4.

$P_0(A|B)$ is the probability of A at $t_0$ on the belief with probability $P_0(B)$ that B will be accepted at $t_1$, whereas in $P_1(A|B)$ B has been accepted at $t_1$. It follows that

$$P_0[P_1(A|B) = P_0(A|B)] = P_0(B).$$

Van Fraassen appears to have mistakenly replaced the '$P_0(B)$' on the right-hand side of this equation by '1'.

$P_0[A|P_1(A)=r]$ in which $P_0[P_1(A) = r] = p$ on one hand and $P_0[A|P_1(A)=s]$ in which $P_0[P_1(A) = s] = 1-p$ on the other, are not mutually exclusive for any $0 \leq r \leq 1$, $0 \leq s \leq 1$, and $0 \leq p \leq 1$. Therefore according to R4 $r = s$ if $r \neq s$. This disproves R4.

Moreover, Van Fraassen has run into the constraint of Section 2.4, which shows that there exists no non-vanishing $P_0(A|B)$ with B='$P_1(A)=r$'. It follows that R1-R4 are all invalid.

## 6.3    *A Purported Coherence Proof*

Several 'coherence proofs' of the mathematical principle of reflection by 'Dutch book' schemes have been proposed (Goldstein, 1983; Van Fraassen, 1984; Gaifman, 1986). Van Fraassen's and Gaifman's proofs are highly informal. Van Fraassen's proof, reduced to essentials, augmented where necessary and in a more general setting which better displays the issue, amounts to this:

Let A be an event whose outcome, true or false, is first accepted at $t_2 > t_1 > t_0$, set B=' $P_1(A)=r$ ', $b = P_{01}(\overline{A}/B) = P_0(\overline{A}|B)$, $\delta = b - (1-r) = r - P_0(A|B)$, and assume that $0 < r < 1$ and $P_0(B) > 0$.

At $t_0$ punter $M$ knows $P_0(A)$, $P_0(B)$, $P_0(AB)$, and $r$, and places with bookie $N$ the following three bets:
  Bet ($a$) pays 1 if $\overline{A}B$ comes true at $t_2$;
  bet ($b$) pays $b$ if $\overline{B}$ comes true at $t_1$;
  bet ($c$) pays $\delta$ if B comes true at $t_1$.
If B comes true at $t_1$ bookie $N$ buys back bet (a) at its updated price $(1-r)$.

The bets cost $P_a = P_0(\overline{A}B) = P_0(B)b$; $P_b = P_0(\overline{B})b$; $P_c = P_0(B)\delta$. Therefore punter $M$ pays bookie $N$ at $t_0$ a total of $P_\Sigma = P_a + P_b + P_c = b + P_0(B)\delta$. If $\overline{B}$ comes true at $t_1$ then $M$ loses $P_\Sigma - b = P_0(B)\delta$, and if B comes true at $t_1$ he loses likewise $P_\Sigma - (1-r) - \delta = P_0(B)\delta$. $P_0(B) > 0$ by assumption. Therefore $\delta = 0$, for otherwise, depending on the predetermined sign of $\delta$, either $M$ or $N$ makes a positive gain whatever the outcome (the bet is a 'Dutch book'; the case is 'incoherent'). Since $P_0(B) > 0$ by assumption, $\delta = 0$ if and only if $r = 1 - b$. But $1 - b = P_0[A | P_0(B)]$ and therefore $\delta = 0$ if and only if $P_0[A | P_0(B)] = r$, which proves R4.

The Proof shows that $P_0[A | P_1(A) = r] = r$ if $0 < r < 1$ is a necessary condition, but another necessary condition is ' $P_0[A | P_1(A) = r] = 0$ ' if $r > 0$ (Section 2.2) with which it is not compatible.

The popularity of Dutch-book proofs stems from their mathematically undemanding nature. They can give necessary conditions but cannot give sufficient conditions, and they are no existence proofs. Dutch-book vulnerability disproves hypotheses but examples of Dutch-book invulnerability do not prove them.[14]

According to the mathematical principle of reflection the imposition of ' $P_0[A | P_1(A) = r] = r$ ' is a necessary and sufficient constraint for the consistency of second order probability. According to the constraint of Section 2.4, the inconsistency that $P_0[A | P_1(A) = r]$ is simultaneously 0 and $r > 0$ if $0 < P_0(B) < 1$, cannot be resolved by the imposition of one of ' $P_0[A | P_1(A) = r] = r > 0$ ' and ' $P_0[A | P_1(A) = r] = 0$ ', unless the rules of probability are suspended, for if one of them is imposed then, by the rules, the other contradicts it.

The Proof is otherwise unsound because condition '$P_0(B) > 0$' is necessary for the inference from '$P_0(B)\delta = 0$' that $\delta = 0$, but $P_0(B) = P_0[P_1(A) = r] \equiv 0$ for epistemic A and pre-determined positive $r$ (Section 2.2).

The Proof has failed.

---

[1] Eur. Ing. Professor Emeritus Amos Nathan, D.Sc., FIET.
Private Scholar
amosnathan1@gmail.com

[2] Instead of $E=\{P(A),P(AB),P(B)\}$ etc., my earlier paper on Probability Dynamics (Nathan, 2006) makes use of '&' as an enumerative connective, as in $E=$'$P(A)\&P(AB)\&P(B)$' etc.

[3] $\overline{A}$ is $\neg A$, etc.

[4] It seems that on p. 80 of Bogen (1995) the phrase 'priors (1) and (3) ... are ... 1' should read -- priors (2) and (3) ... are ... 1' --.

[5] See (Nathan, 2006).

[6] $p \in (0,1)$ stands for $0<p<1$.

[7] The sentence is prompted by Gaifman's (1986) quotation from a cartoon in *The New Yorker*, "There is now 60% chance of rain tomorrow, but, there is 70% chance that later this evening the chance of rain tomorrow will be 80%."

If probabilities are specified to within an integer multiple of some $\Delta$, then the probability of the probability later tonight does not have to vanish but tends to nil with $\Delta$. Similarly, my watch will sometimes show the exact time but the probability that it will do so at some predetermined moment is nil.

[8] In our binary case, $\rho(A,B)=0$ if and only if A and B are independent. In the general case, correlation coefficient $\rho(x,y)$ of random variables *x* and *y* is nil if *x* and *y* are independent, but the converse is not true (Feller, 1957; Gnedenko, 1989).

[9] In (Nathan, 2006) I call it a compound credence.

[10] In (Nathan, 2006) I call it a co-evidence.

[11] A paper comprised of an earlier version of this Section was rejected by *Synthese*, on the strength of a review which asserts that in (Nathan, 2006) the author 'proceeds to *shill* for his own theory'… his theory 'amounts to … *slapping* a real number … on a probability distribution' … 'relative degrees … *whatever that may mean*' … '*arrogant assertions*', etc. For a detailed account see – http://review.amosnathan.com/

[12] The prerequisites of this Section are a few simple relations of PD.

[13] He here refers to it as 'a much maligned principle'.

[14] Sufficiency requires invulnerability to all conceivable Dutch books.